\documentclass[bezier,a4paper,leqno,draft]{amsart}
\usepackage[latin1]{inputenc}
\usepackage{amsmath}
\usepackage{latexsym}
\usepackage{amsthm}
\input{psfig.tex}

\newcommand{\bmath}[1]{\mbox{\mathversion{bold}$#1$}}

\newcommand{\Z}{\bmath{Z}}
\newcommand{\R}{\bmath{R}}

   \newtheorem{theorem}{Theorem}
   \newtheorem{proposition}[theorem]{Proposition}
   \newtheorem{corollary}[theorem]{Corollary}
   \newtheorem{lemma}[theorem]{Lemma}
 \theoremstyle{definition}
   \newtheorem{definition}{Definition}

 \theoremstyle{remark}
   
   \newtheorem*{Rmk}{Remark}
   
\begin{document}

\centerline{\large \bf The first bifurcation point for Delaunay nodoids}

\medskip
\centerline{Wayne Rossman}

\begin{quote} {\bf Abstract:} 
We give two numerical methods for computing the first bifurcation point for 
Delaunay nodoids.  With regard to methods for constructing constant mean 
curvature surfaces, we conclude that the bifurcation point in the analytic 
method of Mazzeo-Pacard is the 
same as a limiting point encountered in the integrable systems method of 
Dorfmeister-Pedit-Wu.  
\end{quote}

\section{Introduction}

Delaunay surfaces in Euclidean $3$-space 
\[ \R^3 = \{(x_1,x_2,x_3) \, | \, x_j \in \R \} \] are constant mean curvature 
(CMC) surfaces of revolution, and they are translationally 
periodic.  By a rigid motion and homothety of $\R^3$ we may place the Delaunay 
surfaces so that their axis of revolution is the $x_1$-axis and their constant 
mean curvature is $H=1$ (henceforth we assume this).  

We consider 
the profile curve in the half-plane $\{(x_1,0,x_3) \in \R^3 \, | \, x_3 > 0\}$ that 
gets rotated about the $x_1$-axis to trace out a Delaunay surface.  This curve 
alternates periodically between maximal and minimal heights (with respect to the 
positive $x_3$ direction), which we refer to as the bulge radius and the 
neck radius, respectively, of the Delaunay surface.  Let us denote the 
neck radius by $r$.  

Delaunay surfaces come in two $1$-parameter families: one is a family of 
embedded surfaces called {\em unduloids} that can be parametrized by the neck radius 
$r \in (0,1/2]$; the other is a family of 
nonembedded surfaces called {\em nodoids} that can be parametrized by the neck radius 
$r \in (0,\infty)$.  For unduloids, $r=1/2$ gives the round cylinder.  For both 
unduloids and nodoids, the limiting singular surface as $r \to 0$ is a chain 
of tangent spheres of radii $1$ centered along the $x_1$-axis.  

In this paper we shall be concerned with nodoids.  We will see that a common 
bifurcation point for Delaunay nodoids is encountered in the following two distinctly 
different approaches for constructing CMC surfaces: 
\begin{enumerate}
\item Using analytic techniques, Mazzeo and Pacard \cite{MP} 
showed existence of a finite value $r_0$ so that for neck radii $r<r_0$ 
the nodoids are nonbifurcating, and for $r>r_0$ the nodoids can bifurcate.  
Bifurcating nodoids are of interest because they deform smoothly through 
families of CMC surfaces 
that are of bounded distance from a fixed line yet are not surfaces of revolution.  
Before the work \cite{MP}, such examples were unknown.  

To study this bifurcation, 
a particular Jacobi operator associated to the second variation formula for 
Delaunay surfaces is used, along with a function space that 
contains only functions with the same translational periodicity as the Delaunay 
surfaces themselves.  Because these functions are translationally periodic, 
they do not have finite $L^2$ norms on the entirety of the surfaces.  Hence 
bifurcation is a different notion from that of "nondegeneracy" of CMC surfaces, 
i.e. CMC surfaces with no nonzero Jacobi fields of finite $L^2$ norm on the 
entire surfaces.  However, the two notions have the common trait of being highly 
useful tools for producing previously unknown CMC surfaces.  
(There has been much work done related to the nondegeneracy 
of CMC surfaces and construction of new CMC examples, see for example the works 
of Kapouleas, Kusner, Mazzeo, Pacard, Pollack, Ratzkin \cite{Kap1}, \cite{Kap2}, 
\cite{KMP}, \cite{MP2}, \cite{MPP}, \cite{JR}.)  As our interest here is in the 
bifurcation of Delaunay nodoids, from the outset we consider only periodic 
functions on Delaunay surfaces.  


Mazzeo and Pacard gave a clear reason for the existence of this bifurcation 
point $r_0$, in terms of the existence of nontrivial nullity for the 
Jacobi operator, but they did not compute the precise value of $r_0$.  
\item Using integrable systems techniques developed by Dorfmeister, Pedit 
and Wu in \cite{dpw}, Dorfmeister, Wu 
\cite{dw} and Schmitt \cite{sch} (see also \cite{kkrs1}) constructed genus $0$ 
CMC surfaces with three asymptotically Delaunay ends.  In \cite{dw} the construction 
was restricted to surfaces with asymptotically unduloidal ends, because such ends 
are embedded.  However, the construction in \cite{sch} and \cite{kkrs1} includes 
asymptotically nodoidal ends as well.  The construction begins with the selection of 
a  certain "DPW potential", and the DPW potential in \cite{sch} and \cite{kkrs1} 
fails to exist when some nodoidal end has an asymptotic neck 
radius greater than $1/2$.  Furthermore, the arguments in \cite{sch} and \cite{kkrs1} 
showing that each end converges to a Delaunay surface work only when the limiting 
Delaunay surface is an unduloid, or is a nodoid with neck size less than $1/2$.  
This suggests that there is possibly some obstruction in the DPW approach to 
CMC surfaces with Delaunay ends that occurs only for asymptotically nodoidal ends 
with asymptotic neck radii at least $1/2$.  

In the construction in \cite{sch} and \cite{kkrs1}, this limiting value $1/2$ 
appears explicitly, but the underlying reasons for its appearance are left 
unexplained.  
\end{enumerate}

It is natural to ask if the bifurcation radius $r_0$ in the first 
approach coincides with the limiting value $1/2$ in the second approach.  The purpose 
of this article is to numerically confirm this:

\begin{quote}
\fbox{
{\bf Numerical result.} The bifurcation radius $r_0$ is equal to $1/2$.
}\end{quote}

This result is of interested with respect to both approaches above.  For the 
first approach, it gives the exact (previously unknown) value for $r_0$.  For the 
second approach, it provides a reason (via the bifurcation properties shown 
by Mazzeo and Pacard) for the existence of the previously mysterious limiting 
value $1/2$.  

In fact, we shall show a stronger numerical result about the first eigenvalue 
of a particular operator, for which the above numerical result is an 
immediate corollary.  

To provide added confidence in the accuracy of our numerical arguments, we give 
two different independent algorithms for showing that $r_0=1/2$ (and for 
showing the stronger numerical result as well).  The first method involves using 
an ordinary differential equation (ODE) solver and symmetry properties of the 
first eigenfunction.  
The second method requires more machinery, the basic tool being the Rayleigh 
quotient characterization for eigenvalues of a self-adjoint operator.  It 
involves numerical integration of smooth bounded functions of a single variable 
on a finite interval, and it 
has the advantage of also giving estimates for other eigenvalues 
beyond the first one.  

\section{Jacobi elliptic functions}

Before introducing equations for nodoids in the next section, we first 
briefly review the properties of Jacobi elliptic functions needed here.  
For $\phi,k \in \R$ with $0 < k^2 < 1$, we set 
\begin{equation}\label{jacobi0} 
\xi = \int_0^\phi \frac{d\psi}{\sqrt{1-k^2 \sin^2 \psi}} \; , 
\end{equation} and then we define 
\[ \text{sn}_{k}(\xi) = \sin \phi \; , \;\;\; 
   \text{cn}_{k}(\xi) = \cos \phi \; , \;\;\; 
   \text{dn}_{k}(\xi) = \sqrt{1-k^2 \sin^2 \phi} \; . \]  Extending these functions 
analytically to $\xi \in \mathbb{C}$, they are defined on the complex plane (with 
singularities).  A short computation gives 
\begin{equation}\label{jacobi1} \frac{d}{d\xi} \text{sn}_k(\xi) = 
\text{cn}_k(\xi) \cdot \text{dn}_k(\xi) \; . \end{equation}  
Defining 
\begin{equation}\label{jacobi2} z = \tilde c \cdot \text{sn}_k(\tilde a \xi+\tilde b) 
\end{equation} for constants $\tilde a, \tilde b$ and 
$\tilde c$, Equation \eqref{jacobi1} then gives 
\begin{equation}\label{jacobi3} \frac{1}{z} \frac{dz}{d\xi} = 
\frac{\tilde a \cdot \text{cn}_k(\tilde a \xi+\tilde b) \cdot 
\text{dn}_k(\tilde a \xi+\tilde b)}{\text{sn}_k(\tilde a \xi+\tilde b)} 
\end{equation} and 
\begin{equation}\label{jacobi4} \left( \frac{dz}{d\xi} \right)^2 + 
\tilde a^2 (1+k)^2 z^2 - \tilde a^2 (\tilde c + k \tilde c^{-1} z^2)^2 = 0 \; . 
\end{equation}
Let $\text{sn}^{-1}_k$ denote the (multi-valued) inverse function of 
$\text{sn}_k$.  Although multi-valued, from here on out 
let us fix $\text{sn}^{-1}_k(1)$ to be equal to the specific value of $\xi \in 
\R$ given by the integral \eqref{jacobi0} with $\phi = \pi/2$.  

\begin{lemma}\label{jacobi5}
Suppose that $\tilde b =${\em sn}$^{-1}_k(1)$, and that 
$\tilde c > 0$ and $\tilde a \in i \R$.  Then the function 
$z$ as in Equation \eqref{jacobi2} is positive (i.e. $z \in \R^+$) and 
periodic with respect to the variable $\xi \in \R$.  
\end{lemma}

\begin{proof}
Using the identity 
\[ \text{sn}_k(\tilde u + \tilde v) = 
\frac{\text{sn}_k(\tilde u) \text{cn}_k(\tilde v) \text{dn}_k(\tilde v) + 
      \text{sn}_k(\tilde v) \text{cn}_k(\tilde u) \text{dn}_k(\tilde u)}
     {1-k^2 \text{sn}_k^2(\tilde u) \text{sn}_k^2(\tilde v)} \] with 
$\tilde u = \tilde a \xi$ and $\tilde v = \tilde b$, we have 
\[ z = \tilde c 
\frac{\text{cn}_k(\tilde a \xi)}{\text{dn}_k(\tilde a \xi)} \; . \]  Then using 
the relations 
\[ \text{cn}_k(i \tilde u) = \frac{1}{\text{cn}_{k^\prime}(\tilde u)} \; , \;\;\; 
\text{dn}_k(i \tilde u) = 
\frac{\text{dn}_{k^\prime}(\tilde u)}{\text{cn}_{k^\prime}(\tilde u)} \] 
with $k^\prime \in \R$ satisfying $(k^\prime)^2+k^2=1$, we have 
\[ z = \frac{\tilde c}{\text{dn}_{k^\prime}(-i \tilde a \xi)} \; . \]  It follows that 
$z$ is positive and periodic when $\xi \in \R$, because $-i \tilde a \xi \in 
\R$ and $\text{dn}_{k^\prime}$ is positive and periodic with respect to a real 
variable.  
\end{proof}

\section{Parametrizing nodoids}

Let $(x_1,x_2,x_3)$ be the usual rectangular coordinates for $\R^3$.  
Consider a Delaunay nodoid with the $x_1$-axis as its axis and with 
constant mean curvature $H=1$.  Let 
\[ (x(t),z(t)) \; , \;\;\; t \in \R \] be a parametrization of the 
profile curve of the nodoid in the $x_1x_3$-plane, and so the surface 
can now be parametrized by \[ {\mathcal D}(t,\theta) = 
(x(t),z(t) \cos \theta, z(t) \sin \theta) \; , 
\;\;\; t \in \R \; , \;\; \theta \in [0,2 \pi ) \; . \] 
Suppose further that the parameter $t$ is chosen to make the 
mapping ${\mathcal D}(t,\theta)$ conformal with respect to the 
coordinates $(t,\theta)$.  Let $t=a$ and $t=b$ be values at which 
the nodoid achieves two adjacent necks, i.e. $z(t)$ has local minima 
at both $t=a$ and $t=b$ equal to the neck radius $r$ 
but at no $t \in (a,b)$.  Conformality implies 
that the first fundamental form is 
\[ ds^2 = ((x^\prime)^2+(z^\prime)^2) dt^2 + z^2 d\theta^2 = 
\rho^2 (dt^2+d\theta^2) \; , \] with 
\[ \rho^2 = (x^\prime)^2+(z^\prime)^2 = z^2 \; . \]  
The second fundamental form is then 
\[ \tfrac{1}{z} (x^{\prime\prime} z^\prime - z^{\prime\prime} x^\prime) dt^2 + 
x^\prime d\theta^2 \; , \] and so the coordinates $(t,\theta)$ 
are curvature line coordinates, that is, the coordinates are isothermic.  
Furthermore, the mean curvature $H=1$ implies 
\[ 2 z^3 - z^\prime x^{\prime\prime} + x^\prime z^{\prime\prime} - 
z x^\prime = 0 \; . \]  This has a first integral, that is, using 
$z^2=(x^\prime)^2 + (z^\prime)^2$ and $z z^\prime = x^\prime x^{\prime\prime} 
+z^\prime z^{\prime\prime}$, 
it is equivalent to \[ (4 x^\prime - 4 z^2)^\prime = 0 \; . \]  Thus 
\begin{equation}\label{weight} m = 4 x^\prime - 4 z^2 \end{equation} 
is satisfied for some $t$-independent constant 
$m$.  Because there are points on a nodoid where $x^\prime = 0$, we 
have \[ m < 0 \; . \]  
In fact, $m \in (-\infty,0)$ is a parameter that determines the full family 
of Delaunay nodoids.  Equation \eqref{weight} can also be established using 
a homology invariant on CMC surfaces called the weight, as explained by 
Korevaar, Kusner and Solomon in 
\cite{KKS}.  We define $m$ as the mass: 

\begin{definition}
Given a nodoid $\mathcal{D}(t,\theta)$ parametrized 
as above satisfying Equation \eqref{weight}, 
we say that $m$ is the {\em mass} (also sometimes called 
the {\em weight} or {\em flux}) of the nodoid $\mathcal{D}(t,\theta)$.  
\end{definition}

Equation \eqref{weight} implies 
\begin{equation}\label{rooteqn} 
z^\prime = \pm \sqrt{z^2-(\tfrac{m}{4}+z^2)^2} \; , \;\;\; m<0 \; . \end{equation}  
Then $x^\prime$ is determined by $z$ via Equation \eqref{weight}, so $x$ is 
determined up to translation, as expected.  

It follows from Equation \eqref{jacobi4} that Equation \eqref{rooteqn} has the 
solution \[ z = -2 B \cdot \text{sn}_{B/A} \left( \chi(t) \right) \; , \;\;\; 
\chi(t) = 2 i A ( t-b ) + \text{sn}^{-1}_{B/A}(1) \; , \] where 
$B=\tfrac{-1}{4}(\sqrt{1-m}-1)$ and $A=\tfrac{1}{2}-B$.  
Lemma \ref{jacobi5} shows us that this solution is positive and periodic with 
respect to $t$, so it is indeed the height function for the profile curve of a 
nodoid in the upper half of the $x_1x_3$-plane.  

From this we can see that the minimum (neck radius) of $z$ is 
\begin{equation}\label{neckradius} r=\sqrt{\frac{1-\tfrac{m}{2}-\sqrt{1-m}}{2}} 
\end{equation} at 
$t=a$ and $t=b$, the maximum (bulge radius) of $z$ is 
\[ \sqrt{\frac{1-\tfrac{m}{2}+\sqrt{1-m}}{2}} \] at 
$t=\tfrac{a+b}{2}$, and $x^\prime=0$ and $z=\sqrt{\tfrac{-m}{4}}$ at both 
$t=\tfrac{3 a+b}{4}$ and $t=\tfrac{a+3 b}{4}$.  

The Gaussian curvature $K$ is determined by $ds^2$ as 
\[ K = \frac{-1}{\rho^2} \triangle (\log \rho) \; , \] 
where \begin{equation}\label{EuclLaplac} 
\triangle = \tfrac{\partial^2}{\partial t^2} + 
\tfrac{\partial^2}{\partial \theta^2} \end{equation} is the standard 
Euclidean Laplacian operator.  

We now introduce a 
function $V$ that will be used in the second variational formula later.  
We set \[ V=(4-2 K) \rho^2 \; . \]  It can be computed that 
\begin{equation}\label{Veqn} V = z^{-2} (2 z^4 + \tfrac{1}{8} m^2) \; . \end{equation} 
Note that $V=V(t)$ is a function of $t$, and is independent of $\theta$.  

\begin{lemma}\label{LemmaOnV}
The function $V=V(t)$ has the following symmetries: \[ V(\ell+t)=V(\ell-t) \;\;\;\; 
\forall \ell \in \{ a, \tfrac{3a+b}{4}, \tfrac{a+b}{2}, \tfrac{a+3b}{4}, b \} \; . \]  
Furthermore, \[ -m \leq V \leq 2-m \;\;\;\; \forall t \in \R \; . \]
\end{lemma}

\begin{figure}[tbp]
\centerline{
        \hbox{
		\psfig{figure=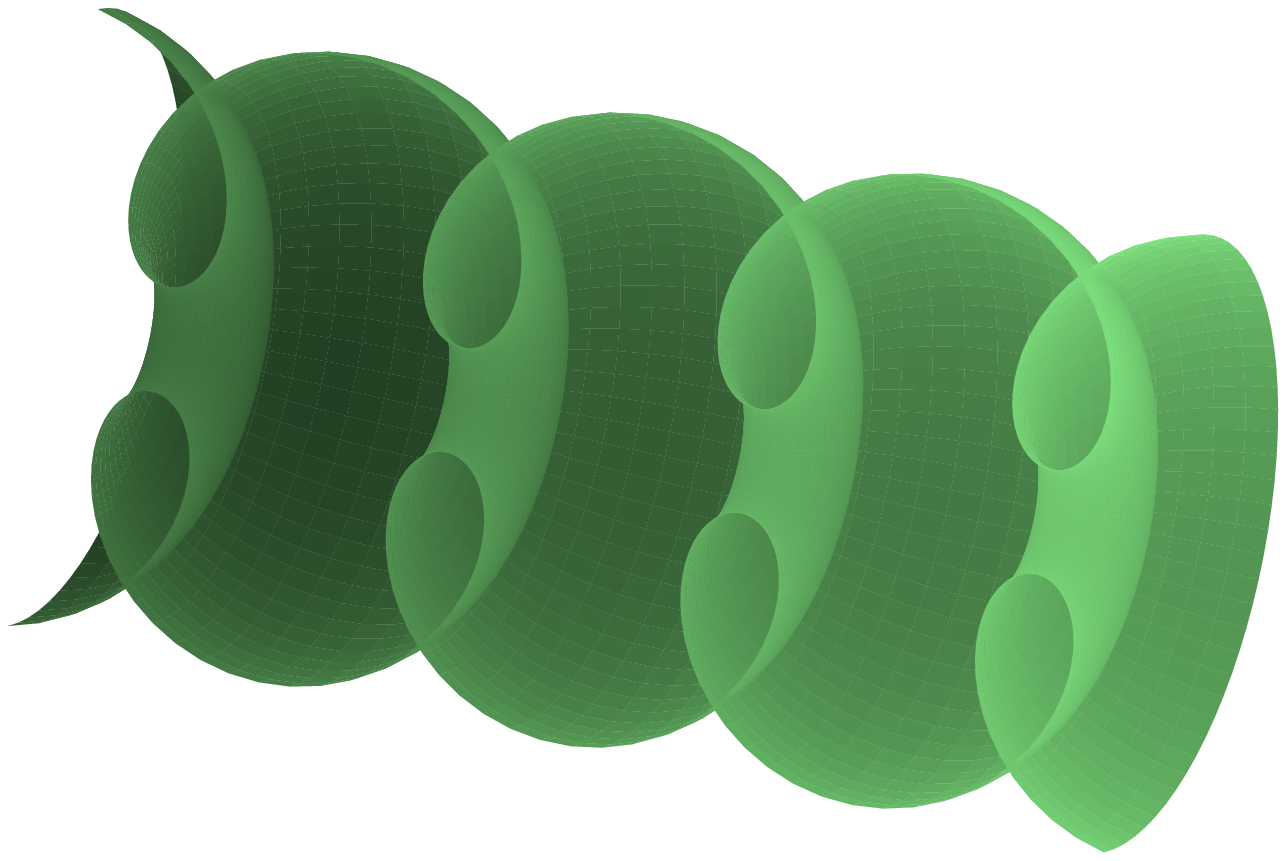,width=2.2in}
		\psfig{figure=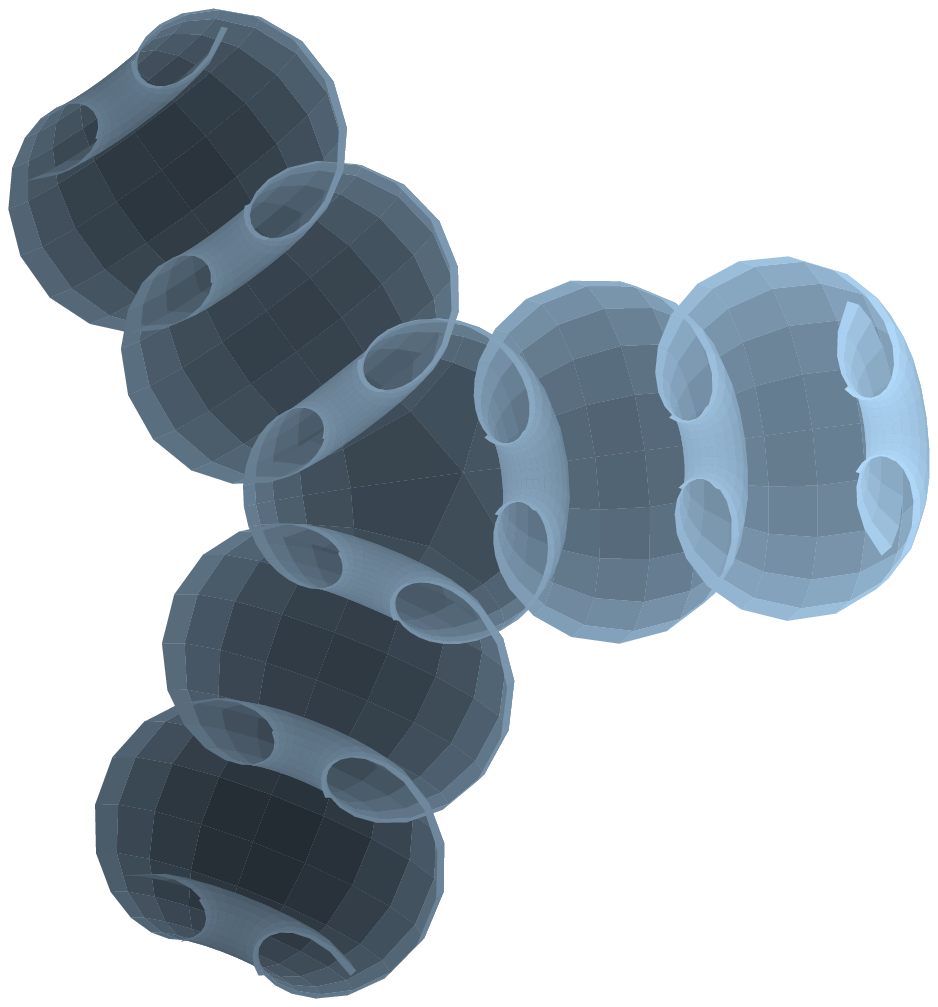,width=2.2in}
	}
}
\caption{\protect\small A nodoid and a CMC surface with three ends asymptotic 
to nodoids.  One of two halves of each surface is shown here.  
(The computer graphics were made using N. Schmitt's cmclab software 
\cite{cmclab} and K. Polthier's Javaview software \cite{P1}.)}
\label{fig:1}
\end{figure}

\begin{proof}
We will show that $V(a+t)=V(a-t)$ and 
$V(\tfrac{3a+b}{4}+t)=V(\tfrac{3a+b}{4}-t)$, and then all other symmetries in 
the lemma follow.  Because $z(t)$ itself has the symmetry 
$z(a+t)=z(a-t)$, clearly also $V(a+t)=V(a-t)$ by Equation \eqref{Veqn}.  
To show $V(\tfrac{3a+b}{4}+t)=V(\tfrac{3a+b}{4}-t)$, 
we first note, using Equation \eqref{jacobi3}, that 
\[ \frac{z^\prime(t)}{z(t)} = 2 i A 
\frac{\text{cn}_{B/A} \left( \chi(t) \right) \text{dn}_{B/A} \left( \chi(t) \right)}
     {\text{sn}_{B/A} \left( \chi(t) \right)} \; , \] 
and so \begin{equation}\label{insideproof} 
\frac{z^\prime(\tfrac{3a+b}{4}+t)}{z(\tfrac{3a+b}{4}+t)}= 
\frac{z^\prime(\tfrac{3a+b}{4}-t)}{z(\tfrac{3a+b}{4}-t)} \; . \end{equation} 
Substituting $m^2/8 = (2-m) z^2 -2 z^4 -2 (z^\prime)^2$ from Equation 
\eqref{rooteqn} into Equation \eqref{Veqn}, we have 
\begin{equation}\label{Veqninsideproof} 
V(t) = 2-m-2 \left( \tfrac{z^\prime}{z} \right)^2 \; . \end{equation} 
Then by Equation \eqref{insideproof}, we have 
$V(\tfrac{3a+b}{4}+t)=V(\tfrac{3a+b}{4}-t)$.  
Equation \eqref{Veqninsideproof} also implies that $V \leq 2-m$ for all $t$.  

To see that $V \geq -m$ for all $t$, we simply consider $V=2 \zeta + \tfrac{m^2}{8} 
\tfrac{1}{\zeta}$ as a function of $\zeta = z^2 > 0$.  Elementary calculus gives that 
$2 \zeta + \tfrac{m^2}{8 \zeta} \geq -m$ for all $\zeta > 0$.  
\end{proof}

\section{Second variation for nodoids}

We now consider an arbitrary volume-preserving periodic smooth variation 
of the surface, that is, a family of surfaces ${\mathcal D}_s(t,\theta)$ 
that is 
\begin{enumerate}
\item smooth with respect to $s \in (-\epsilon, \epsilon)$ for some 
$\epsilon > 0$, 
\item an immersion with respect to the coordinates $(t,\theta)$ for any fixed 
$s \in (-\epsilon, \epsilon)$, 
\item satisfying 
${\mathcal D}_0(t,\theta) = {\mathcal D}(t,\theta)$ for all $t \in \R$ 
and $\theta \in [0,2\pi)$, 
\item periodic with the same period for all $s \in (-\epsilon, \epsilon)$, 
i.e. ${\mathcal D}_s(t_1,\theta)={\mathcal D}_s(t_2,\theta)$ for all 
$s \in (-\epsilon, \epsilon)$ and all $t_2-t_1$ an integer multiple of $b-a$, 
\item with compact regions bounded by the surfaces ${\mathcal D}_s(t,\theta), t \in 
[a,b], \theta \in [0,2\pi)$ and the two disks 
\[ D_a(s) = \{(x(a),x_2,x_3) \, | \, x_2^2+x_3^2 \leq (z(a)+\mathcal{O}(s))^2 \} \; , \] 
\[ D_b(s) = \{(x(b),x_2,x_3) \, | \, x_2^2+x_3^2 \leq (z(b)+\mathcal{O}(s))^2 \} \] 
having the same volume for all $s \in (-\epsilon,\epsilon)$.  
\end{enumerate}
The disks $D_a(s)$ and $D_b(s)$ are allowed to deform smoothly in $s$, and are not 
necessarily perfectly round when $s \neq 0$.  Thus these deformations do not have 
any fixed "boundary" curves.  The essential properties of these deformations are only 
that they are both periodic and volume-preserving.  

Let $\vec{N}$ denote a unit normal vector to ${\mathcal D}_0(t,\theta)$.  
We define the function $u$ by 
\[ u=u(t,\theta)= \left\langle \left. 
\tfrac{d}{ds} {\mathcal D}_s(t,\theta) \right|_{s=0} , \vec{N} \right\rangle 
\in \mathcal{F} \; . \] 
The volume-preserving condition 5. above implies 
\[ \int_a^b \int_0^{2\pi} u \rho^2 d\theta dt = 0 \; . \]  

\begin{definition}
We will call the compact portion ${\mathcal D}(t,\theta), t \in 
[a,b], \theta \in [0,2\pi)$ a {\em fundamental piece} of the nodoid.  
\end{definition}

Let $A(s)$ denote the area of the surface ${\mathcal D}_s(t,\theta)$ for $t \in 
[a,b], \theta \in [0,2\pi)$.  Because the nodoid is CMC $H=1$, the first variation 
formula for a fundamental piece is 
\[ \left. \frac{d}{ds} A(s) \right|_{s=0} = 
\int_a^b \int_0^{2\pi} u \rho^2 d\theta dt = 0 \; . \]  
Thus it is the second variation formula for volume-preserving variations 
(see \cite{bc}) 
\[ \left. \frac{d^2}{ds^2} A(s) \right|_{s=0} = 
\int_a^b \int_0^{2\pi} u (-\triangle_{ds^2} u - (4 - 2 K) u) \rho^2 d\theta dt \] 
\[ = \int_a^b \int_0^{2\pi} u \cdot {\mathcal L}(u) d\theta dt = 0 \; , \;\;\;\;\; 
{\mathcal L}(u) := -\triangle u - V u \; , \] 
with $\triangle_{ds^2}$ (respectively $\triangle$) the Laplace-Beltrami 
operator determined by $ds^2$ (respectively the Euclidean Laplacian as in 
\eqref{EuclLaplac}), that will determine if the variation increases or decreases area 
(when this formula is nonzero).  We are applying the first and second variation 
formulas only to periodic variations here, but we remark that they 
hold for other nonperiodic types of variations as well.  

\section{Spherical harmonics}

We now consider the eigenvalue problem for $\mathcal L$.  
We first introduce the function space 
\[ \mathcal{F} = \{ u=u(t,\theta) \in C^\infty(t,\theta) \, | \, 
      u(t_1,\theta)=u(t_2,\theta) \text{ for }t_1-t_2 \in (b-a) \cdot \Z \} . \]
The eigenvalue problem is to find $u \in \mathcal{F}$ and $\lambda \in \R$ such that 
\begin{equation}\label{eigfcteqn} 
{\mathcal L}(u) = \lambda u \; . \end{equation} 
Let us decompose such an eigenfunction $u$ into its spherical harmonics: 
$u=u(t,\theta)$ can be written as 
\[ u = u_0(t) + \sum_{j \geq 1} (u_{j,+}(t) \cdot \cos(j\theta) + 
u_{j,-}(t) \cdot \sin(j\theta)) \; , \] where $u_0(t)$, $u_{j,+}(t)$ and 
$u_{j,-}(t)$ are periodic functions of $t$, that is, they lie in 
the smaller function space 
\[ \hat{\mathcal F} = 
\mathcal{F} \cap \{ u=u(t,\theta) | \tfrac{\partial}{\partial \theta} 
u \equiv 0 \} \; , \]  
i.e. those functions in $\mathcal F$ that are independent of $\theta$.  
Defining the operators 
\[ {\mathcal L}_j = -\tfrac{\partial^2}{\partial t^2} - V + j^2 \] for 
$j \in \Z^+ \cup \{ 0 \}$ on the function space 
$\hat{\mathcal F}$, the relation \eqref{eigfcteqn} gives also that 
\[ {\mathcal L}_0(u_0) = \lambda u_0 \; , \;\; {\mathcal L}_j(u_{j,\pm}) = 
\lambda u_{j,\pm} \] for all $j \geq 1$.  So the following lemma is immediate: 

\begin{lemma}
A real number $\lambda$ is an eigenvalue of $\mathcal{L}$, i.e. 
$\mathcal{L}(u)=\lambda u$ for some $u \in \mathcal{F}$, if and only if 
$\lambda$ is also an eigenvalue of $\mathcal{L}_j$ for at least one value 
of $j \in \Z^+ \cup \{ 0 \}$, i.e. $\mathcal{L}_j(\hat u)=\lambda \hat u$ 
for some $\hat u \in \hat{\mathcal F}$ and some $j \in \Z^+ \cup \{ 0 \}$.  
\end{lemma}

The following lemma is analogous to Proposition 4.4 in \cite{MP}, but different 
notations were used there:  

\begin{proposition}\label{zeroandminusone}
Both $-1$ and $0$ are eigenvalues of $\mathcal{L}_0$.  
\end{proposition}

\begin{proof}
We first show that $0$ is an eigenvalue.  
Considering the vector $(1,0,0)$ as a 
constant vector field of $\R^3$, the associated 
translational flow gives a periodic 
volume-preserving deformation $\mathcal{D}_s$ of the nodoid 
$\mathcal{D}=\mathcal{D}_0$.  Because this flow is actually a family of 
rigid motions, it is also 
area-preserving.  It is a classical fact that the scalar product of the 
unit normal vector of a CMC surface with a Killing field is in the null-space 
of the Jacobi operator of the surface (see \cite{Choe}, page 196, or the proof
of Theorem 2.7 in \cite{bgs}).  Hence, if $(1,0,0)$ is decomposed into 
$(1,0,0)=u \vec{N} + \vec{v}$ at each point of $\mathcal{D}$ with $u \vec{N}$ 
normal to $\mathcal{D}$ and $\vec{v}$ tangent to $\mathcal{D}$, and with 
$u \in \hat{\mathcal F}$ and $|\vec{N}|=1$, then $\mathcal{L}(u)=0$.  
Since $u \in \hat{\mathcal F}$, it follows that $0$ is an eigenvalue of 
$\mathcal{L}_0$ for any choice of $m < 0$.  Also, by direct computation, we 
have $u = -z^\prime/z$, and then using $(z^\prime)^2=z^2 - (z^2+m/4)^2$ and 
$z^{\prime\prime} = z-2 z (z^2+m/4)$, we have $\mathcal{L}_0(u)=0$.  

Now we show that $-1$ is an eigenvalue.  
Similarly to the previous case with eigenvalue $0$, we now consider the 
vector $(0,0,1)$, 
producing a smooth vector field of $\R^3$, whose associated 
translational flow again gives a periodic 
volume-preserving deformation $\mathcal{D}_s$ of $\mathcal{D}=\mathcal{D}_0$ 
that is once again area-preserving.  Now we decompose 
$(0,0,1)=u \sin \theta \vec{N} + \vec{v}$ at each point of $\mathcal{D}$ 
with $u \sin \theta \vec{N}$ 
normal to $\mathcal{D}$ and $\vec{v}$ tangent to $\mathcal{D}$. Then 
$\mathcal{L}(u \sin \theta)=0$, and $u 
\in \hat{\mathcal F}$ is an eigenfunction of $\mathcal{L}_1$ with eigenvalue 
$0$.  It follows that $-1$ is an eigenvalue of $\mathcal{L}_0$ for any choice 
of $m < 0$.  Again, we could also see this by direct computation: we 
have $u = x^\prime/z=z+m/(4z)$, and then using $(z^\prime)^2=z^2 - (z^2+m/4)^2$ and 
$z^{\prime\prime} = z-2 z (z^2+m/4)$, we have $\mathcal{L}_0(u)=-u$.  
\end{proof}

Both the operator $\mathcal L$ on the function space $\mathcal F$ and 
the operators $\mathcal{L}_j$ on the function space $\hat{\mathcal F}$ 
are essentially self-adjoint, and hence standard functional analysis 
arguments (see \cite{be} or \cite{u} or \cite{chavel} 
for example) give the following result: 

\begin{proposition}\label{eigvalproperties}
Each of the operators $\mathcal L$ on the function space $\mathcal F$ and 
$\mathcal{L}_j$ on the function space $\hat{\mathcal F}$ satisfy the following: 
\begin{enumerate}
\item all eigenvalues are real, and there are a countably infinite number of them, 
\item all eigenvalues are greater than some real constant, 
\item the eigenvalues do not accumulate at any finite real value, 
\item when written in increasing order, the eigenvalues increase to $+\infty$, 
\item the eigenspace associated to each eigenvalue is finite dimensional, 
\item the collection of eigenspaces spans the full function space.  
\end{enumerate}
\end{proposition}

We refer to the eigenvalue that is less than all the others as the first eigenvalue.  
The second eigenvalue is the one that is less than all but the first eigenvalue.  
The third eigenvalue is the one that is less than all but the first and 
second eigenvalues, and so on.  

The Courant nodal domain theorem, which is valid in our setting (see \cite{cheng}), 
tells us that the number of nodal domains of the eigenfunction associated to 
the $k$'th eigenvalue is at most $k$ (here we are counting eigenvalues with 
multiplicity).  If the eigenspace associated to the first 
eigenvalue contained two linearly 
independent functions, then some linear combination of the 
two functions would be a function that attains both positive 
and negative values.  This would contradict the Courant nodal domain theorem, so the 
eigenspace associated to the first eigenvalue is $1$ dimensional.  

Also, in the case of the operator $\mathcal{L}_j$, the domain on which the functions 
are defined is $1$ dimensional, so $\mathcal{L}_j$ is an ordinary differential 
equation, not a partial differential equation.  
Hence, saying that any 
(not identically zero) eigenfunction associated to the first eigenvalue 
has at most one nodal domain is equivalent to saying that any such eigenfunction 
never attains the value zero.  We conclude the following: 

\begin{proposition}\label{courant}
For both the operator $\mathcal{L}$ defined on the function space $\mathcal{F}$ and the 
operator $\mathcal{L}_j$ defined on the function space $\hat{\mathcal F}$, 
the eigenspace associated to the first eigenvalue is $1$ dimensional.  Furthermore, 
in the case of $\mathcal{L}_j$, 
any not-identically-zero function in this eigenspace never attains the value zero.  
\end{proposition}

\begin{Rmk}
The eigenvalues $-1$ and $0$ of $\mathcal{L}_0$ are actually the second and 
third eigenvalues of $\mathcal{L}_0$.  This was shown in \cite{MP} by counting 
the nodal domains of the corresponding eigenfunctions.  
\end{Rmk}

\section{Primary result}

The first two operators $\mathcal{L}_0$ and $\mathcal{L}_1$ never give 
any bifurcation in the sense of \cite{MP}.  
The first bifurcation point (as in \cite{MP}) is defined as follows: 

\begin{definition}\label{bifptdefn}
The {\em first bifurcation point} 
for Delaunay nodoids occurs at the value of $m < 0$ closest to zero 
for which ${\mathcal L}_j$ for some $j \geq 2$ 
has a nonpositive eigenvalue, i.e. 
\[ \mathcal{L}_j u = \lambda u \;\; \text{with} \;\; \lambda \leq 0 
\;\; \text{and} \;\; u \in \hat{\mathcal F} \;\; \text{for some} \;\; j \geq 2 \; . \]  
\end{definition}

This occurs at the largest value of $m$ 
for which the first eigenvalue of ${\mathcal L}_2$ is zero, and is equivalent 
to the first eigenvalue $\lambda_0$ of ${\mathcal L}_0$ being $-4$.  This 
$\lambda_0$ can be computed using a minimum of the Rayleigh quotient (\cite{be}, 
\cite{u}, \cite{chavel}): 
\begin{equation}\label{Rayleigh} \lambda_0 = \min_{u(t) \in \hat{\mathcal F} 
\setminus \{0\} } 
\frac{\int_a^b u(t) {\mathcal L}_0(u(t)) dt}{\int_a^b u(t)^2 dt} \; . \end{equation}

\begin{Rmk}
Bifurcation points are points where bifurcation actually occurs, that is, in 
any neighborhood of the Delaunay surface there are CMC surfaces that are not 
rotationally symmetric.  
What we have defined in Definition \ref{bifptdefn} is an "infinitesimal symmetry 
breaking point" 
where bifurcation might happen, because having zero as an eigenvalue 
of ${\mathcal L}_2$ is a necessary condition for the bifurcation to occur.  
However, because the multiplicity of the zero eigenvalue of 
${\mathcal L}_2$ is one (see Proposition \ref{courant}), and because the 
derivative of the first eigenvalue with respect to the parameter $m$ is not zero 
(see the primary numerical result below), it turns out that the 
infinitesimal symmetry breaking point is an actual symmetry breaking 
point, i.e. a true bifurcation point \cite{Fpersonal}.  
\end{Rmk}

\begin{figure}[tbp]
\centerline{
        \hbox{
		\psfig{figure=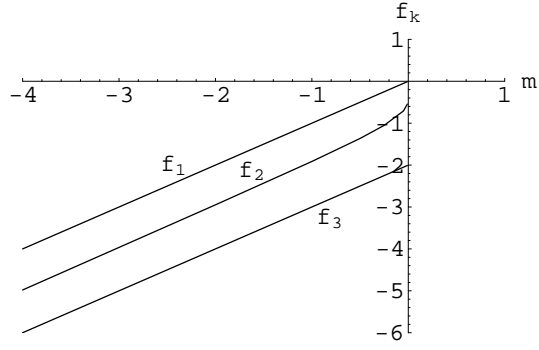,width=3.0in}
	}
}
\caption{\protect\small Plots of the functions $f_1=m$, $f_3=m-2$ and 
$f_2=\frac{-1}{b-a} \int_a^b V dt$.  $\lambda_0$ is also a function of $m$, and 
Lemmas \ref{simplelemma1} and \ref{simplelemma2} imply 
that $\lambda_0$ must lie between $f_2$ and $f_3$.}
\label{fig:2}
\end{figure}

We will give our primary numerical result and two different numerical methods for 
showing it.  But first let us give two simple mathematically rigorous lemmas that 
support this upcoming numerical result.  Both of these results are immediate from 
the Rayleigh quotient formulation for $\lambda_0$ above.  The first lemma is 
also shown in Proposition 4.4 of \cite{MP}, with different notations.  

\begin{lemma}\label{simplelemma1}
For a nodoid of mean curvature $1$ and mass $m$, 
the first eigenvalue $\lambda_0$ of $\mathcal{L}_0$ satisfies 
$m-2 \leq \lambda_0 \leq m$.  
\end{lemma}

\begin{proof}
Using the Rayleigh quotient characterization \eqref{Rayleigh}, this follows 
directly from the fact that $-m \leq V=V(t) \leq 2-m$ for all $t$, as shown in 
Lemma \ref{LemmaOnV}.  
\end{proof}

We then immediately have: 

\begin{corollary}
For nodoids of mean curvature $1$ and mass $m$, the first 
bifurcation point occurs for some $m$ between $-4$ and $-2$.  
\end{corollary}

\begin{lemma}\label{simplelemma2}
For a nodoid of mean curvature $1$ and mass $m$, 
the first eigenvalue $\lambda_0$ of $\mathcal{L}_0$ satisfies 
$\lambda_0 \leq \frac{-1}{b-a} \int_a^b V dt$.  
\end{lemma}

\begin{proof}
Inserting $u \equiv 1$ into the Rayleigh quotient in Equation \eqref{Rayleigh}, this 
lemma follows.  
\end{proof}

Numerically checking that $\frac{-1}{b-a} \int_a^b V dt$ is an increasing 
function of $m<0$ that becomes $-4$ when 
(and only when) $m$ is approximately $-3.036$, we have the following: 

\begin{quote}
{\bf Preliminary numerical result.}  For nodoids of mean curvature $1$ and mass 
$m$, the first bifurcation point occurs for some $m \geq -3.036$.  
\end{quote}

Now, armed with the knowledge that the bifurcation point must occur for some 
$m$ between $-3.036$ and $-2$, we state our primary numerical result: 

\smallskip

\begin{quote}
{\bf Primary numerical result.}  For a nodoid of mean curvature $1$ and mass $m$, 
the first eigenvalue $\lambda_0$ of $\mathcal{L}_0$ is 
\[ \lambda_0=m-1 \; . \]  Hence the first bifurcation 
point for nodoids occurs when $m=-3$.  
\end{quote}

\medskip

By the formula \eqref{neckradius}, $m=-3$ precisely when the 
neck radius is $r=1/2$, thus the numerical result stated in the introduction is a 
direct corollary of this numerical result here.  

\section{First method}

In this and the next section we give two independent methods for 
numerically confirming our primary numerical result.  The first method 
in this section is the simpler of the two.  

We wish to solve
\begin{equation}\label{method1eqn} \frac{d^2}{dt^2} u +
(V+\lambda) u = 0 \end{equation} for some
$\lambda \in [m-2,m]$ and some function $u = u(t) \in \hat{\mathcal F}$.
If there is only one such value for $\lambda<-1$ in the range $[m-2,m]$ (and 
we will see that this is the case), then this $\lambda$ will be the first 
eigenvalue $\lambda_0$.  We begin with a mathematically rigorous lemma that 
this first numerical method is based on: 

\begin{lemma}\label{Vsymmetry}
We fix any $m<0$ and consider the eigenvalue equation \eqref{method1eqn}.  
Let $u$ be a not-identically-zero 
eigenfunction associated to the first eigenvalue $\lambda_0$.  Then 
$u$ cannot attain zero at any value of $t$, and 
\[ u(\ell+t)=u(\ell-t) \;\;\;\; \text{for all} \;\; \ell \in \{ a, 
\tfrac{3a+b}{4}, \tfrac{a+b}{2}, \tfrac{a+3b}{4}, b \} \; . \]  
In particular, \[ u^\prime(a)=u^\prime(\tfrac{3a+b}{4})=u^\prime(\tfrac{a+b}{2})=
u^\prime(\tfrac{a+3b}{4})=u^\prime(b)=0 \; . \]  
\end{lemma}

\begin{proof}
Proposition \ref{courant} tells us that the eigenspace associated 
to the first eigenvalue $\lambda_0$ must be one dimensional.  
If $u$ did not have all of the 
same symmetries as $V$ as in Lemma \ref{LemmaOnV}, then there would be some 
$\ell \in \{ a, 
\tfrac{3a+b}{4}, \tfrac{a+b}{2}, \tfrac{a+3b}{4}, b \}$ so that $u(\ell+t)$ and 
$u(\ell-t)$ are two linearly independent eigenfunctions associated to the 
same first eigenvalue $\lambda_0$.  This is a contradiction, implying the lemma.  
\end{proof}

\begin{figure}
 \begin{center}
  \begin{tabular}{ccc}
		\psfig{figure=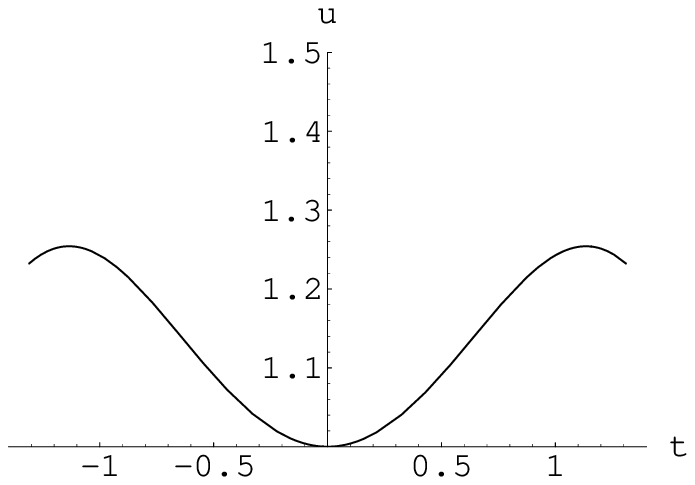,width=1.9in} &
		\psfig{figure=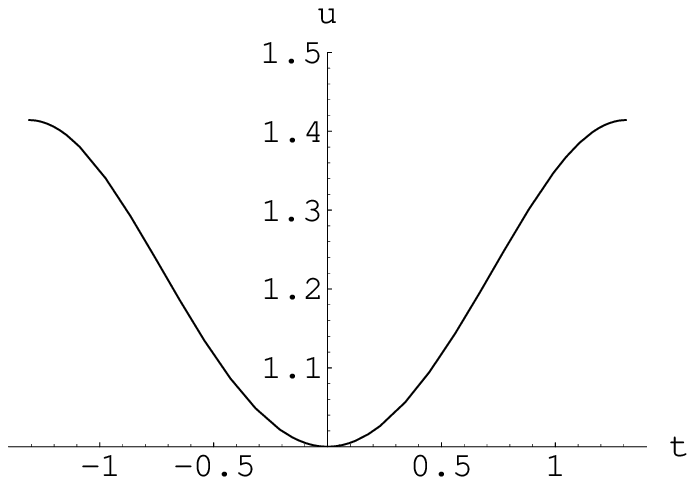,width=1.9in} &
		\psfig{figure=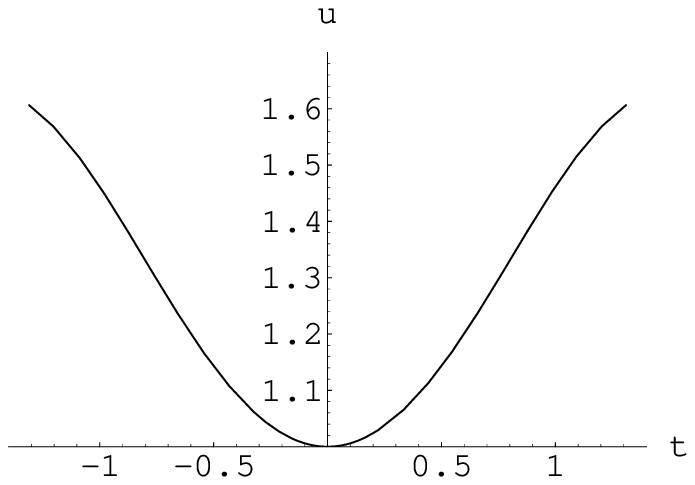,width=1.9in} \\
   \small $m=-1$, $\lambda=m-0.8$ &
   \small $m=-1$, $\lambda=m-1.0$ &
   \small $m=-1$, $\lambda=m-1.2$ 
  \end{tabular}
  \begin{tabular}{ccc}
		\psfig{figure=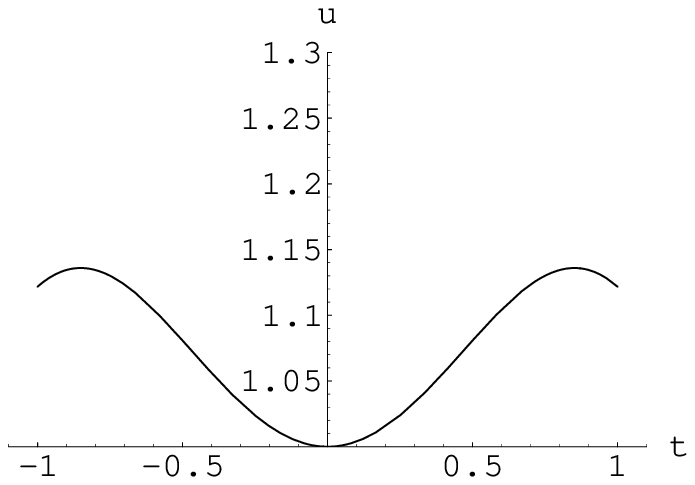,width=1.9in} &
		\psfig{figure=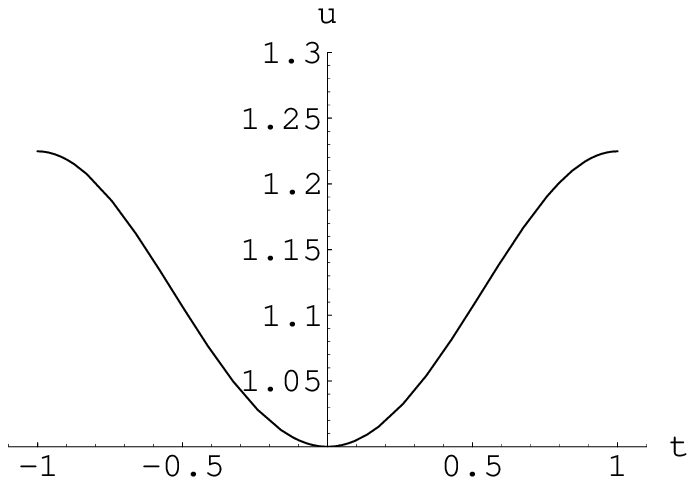,width=1.9in} &
		\psfig{figure=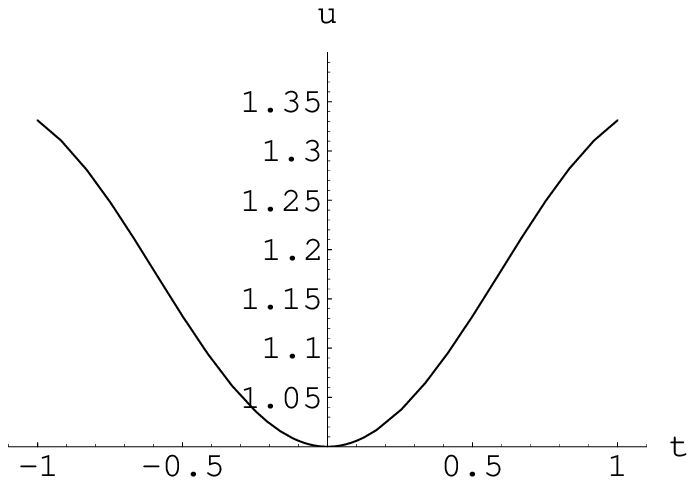,width=1.9in} \\
   \small $m=-2$, $\lambda=m-0.8$ &
   \small $m=-2$, $\lambda=m-1.0$ &
   \small $m=-2$, $\lambda=m-1.2$ 
  \end{tabular}
  \begin{tabular}{ccc}
		\psfig{figure=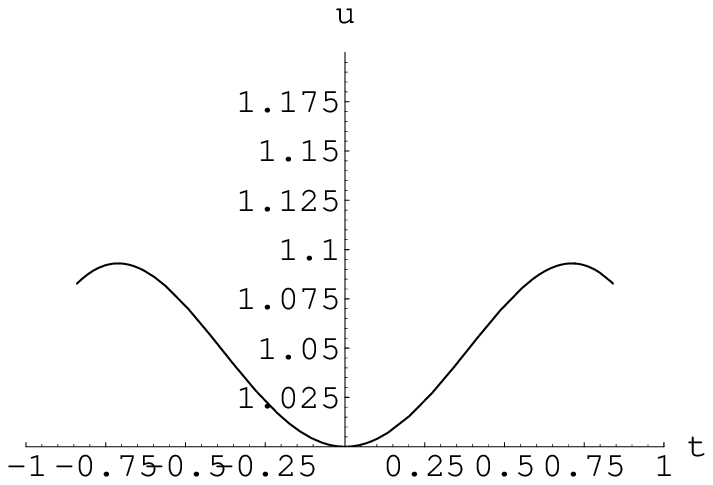,width=1.9in} &
		\psfig{figure=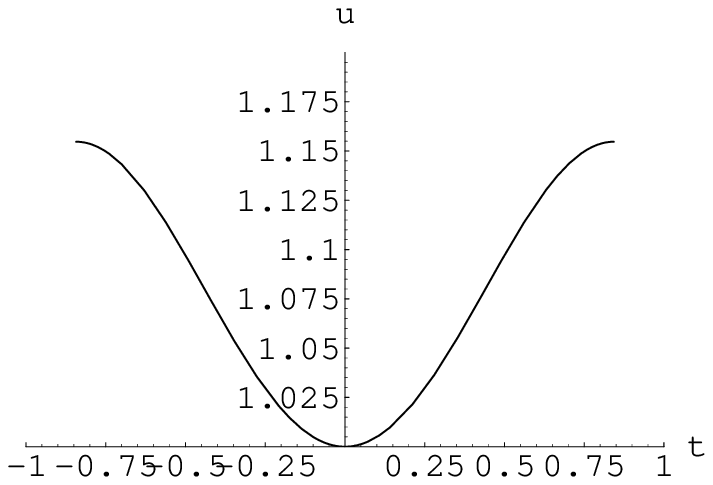,width=1.9in} &
		\psfig{figure=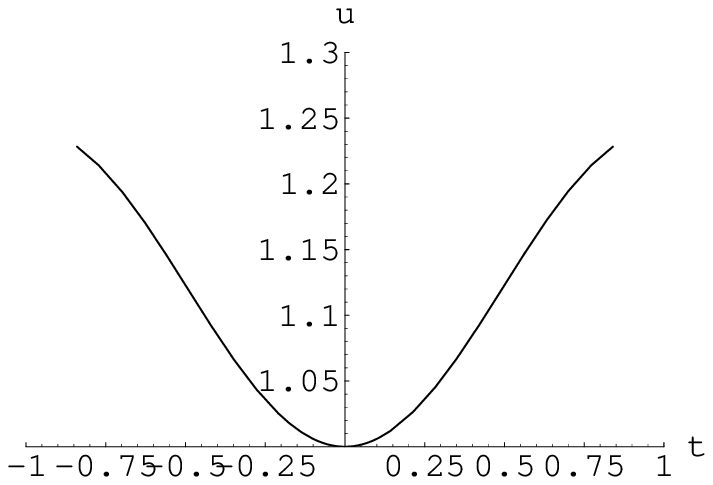,width=1.9in} \\
   \small $m=-3$, $\lambda=m-0.8$ &
   \small $m=-3$, $\lambda=m-1.0$ &
   \small $m=-3$, $\lambda=m-1.2$ 
  \end{tabular}
 \end{center}
\caption{\protect\small The function $u$ computed numerically 
via the first method.  In each case, the graph 
is drawn for $t \in [a,\tfrac{a+b}{2}]$ with $\tfrac{b}{3}=-a=1.3108$ (resp. 
$1.0$, $0.8428$) for $m=-1$ (resp. $-2$, $-3$).}
\label{fig:3}
\end{figure}

Because multiplying $u$ by a
real scalar does not affect Equation \eqref{method1eqn}, we may assume
that \[ u(\tfrac{3 a+b}{4})=1 \; . \]  
The numerical method is then to numerically solve Equation
\eqref{method1eqn} for $u$ with initial conditions
\[ u(\tfrac{3 a+b}{4})=1 \; , \;\;\; u^\prime(\tfrac{3 a+b}{4})=0 
\] by some ODE solver
(such as "NDSolve" in Mathematica), and find the value of
$\lambda$ that gives a periodic solution $u$, that is, that gives 
\[ u(a)=u(\tfrac{a+b}{2})=u(b) \; , \;\;\; u^\prime(a)=u^\prime(\tfrac{a+b}{2})= 
u^\prime(\tfrac{a+3b}{4})=u^\prime(b)=0 \; , \;\;\; u(\tfrac{a+3 b}{4})=1 \; , 
\;\;\; u(t) > 0 \; \forall t \in [a,b] \; , \]
or equivalently, by the symmetries of $V$ in Equation \eqref{method1eqn} and by 
the symmetries of $u$ in Lemma \ref{Vsymmetry}, that simply gives 
\begin{equation}\label{closingcond} 
u^\prime(a)=u^\prime(\tfrac{a+b}{2})= 0 \; . \end{equation} 
Numerically finding this value of $\lambda$ for various values 
of $m<0$, one finds that it is always 
$m-1$, verifying the primary numerical result.  

Figure \ref{fig:3} demonstrates how this occurs, for the cases $m=-1$, $-2$ and $-3$.  
In each case, the ODE solver produces a solution that satisfies 
Equation \eqref{closingcond} 
precisely when $\lambda = m-1$.  When $\lambda \in (m-1,m]$ (for example, 
$\lambda = m-0.8$ as in Figure \ref{fig:3}), 
the values of $t$ ($\neq \tfrac{3a+b}{4}$) closest 
to $\tfrac{3a+b}{4}$ where $u^\prime(t)=0$ satisfy $|t-\tfrac{3a+b}{4}| < 
\tfrac{b-a}{4}$.  When $\lambda \in [m-2,m-1)$ (for example, 
$\lambda = m-1.2$ as in 
Figure \ref{fig:3}), the values of $t$ ($\neq \tfrac{3a+b}{4}$) closest 
to $\tfrac{3a+b}{4}$ where $u^\prime(t)=0$ satisfy $|t-\tfrac{3a+b}{4}| > 
\tfrac{b-a}{4}$.  We need precisely $u^\prime(\tfrac{a+b}{2})=u^\prime(a)=0$, 
and this occurs exactly when $\lambda=m-1$.  This has been checked for numerous 
other values of $m<0$ as well.  

\section{Second method}

The second method we present here is more complicated than the
first one, but we wish to consider it for the following two
reasons: 1) it does not use an ODE solver and the algorithm 
is independent of the first method above, thus it gives a second
fully-independent confirmation of the numerical result; 2) it is a stronger
method in that it also gives estimates for other eigenvalues
of $\mathcal{L}_0$, not just the first one.  

As we saw in Proposition \ref{zeroandminusone} and the remark just after Proposition 
\ref{courant}, in addition to $\lambda_0$,
$-1$ and $0$ are both also eigenvalues of $\mathcal{L}_0$ for all 
$m<0$, and $-1$ and $0$ are 
actually the second and third eigenvalues of $\mathcal{L}_0$.  
The fact that this 
second method precisely estimates the second and third eigenvalues $-1$ and $0$
gives us confidence that the method is accurately estimating $\lambda_0$ as
well.

We choose a basis $\{ \mathcal{B}_j=\mathcal{B}_j(t) \}_{j = 1}^\infty$ for 
$\hat{\mathcal F}$ as 
\[ \mathcal{B}_1 = \frac{1}{\sqrt{b-a}} \; , \] 
\[ \mathcal{B}_j = \sqrt{\frac{2}{b-a}} \cdot \cos \left( \frac{\pi j (t-a)}{b-a} 
    \right) \;\;\; \text{ for even } \; j \in 2 \Z^+ \; , \]  
\[ \mathcal{B}_j = \sqrt{\frac{2}{b-a}} \cdot \sin \left( \frac{\pi (j-1) (t-a)}{b-a} 
    \right) \;\;\; \text{ for odd } \; j \in 2 \Z^+ + 1 \; . \]  
This is an orthonormal basis for $\hat{\mathcal F}$ with respect to the 
Euclidean $L^2$ norm on the interval $[a,b]$.  Any $u \in \hat{\mathcal F}$ 
can be expanded as 
\[ u = \sum_{j=1}^\infty a_j \mathcal{B}_j \] for constants $a_j \in \R$.  
The Rayleigh quotient characterization \eqref{Rayleigh} then gives 
\[ \lambda_0 = \min_{\sum_{j \geq 1} a_j^2 > 0} 
\frac{\sum_{j,k \geq 1} a_j a_k \alpha_{jk}}{\sum_{j,k \geq 1} a_j a_k} \; , \]
where 
\[ \alpha_{jk} = \int_a^b -\mathcal{B}_j \tfrac{d^2}{dt^2}\mathcal{B}_k dt - 
   \int_a^b V \mathcal{B}_j \mathcal{B}_k dt \; . \]  

The integrals $\int_a^b -\mathcal{B}_j \tfrac{d^2}{dt^2}\mathcal{B}_k dt$ can be 
explicitly computed as
\[ \int_a^b -\mathcal{B}_j \tfrac{d^2}{dt^2}\mathcal{B}_k dt = 
\delta_{jk} \left[ \frac{j}{2} \right]^2 \frac{4 \pi^2}{(b-a)^2} \; , \] 
where $\delta_{jk}$ is the Kronecker-delta function and $\left[ \tfrac{j}{2} 
\right]$ is the greatest integer less than or equal to $\tfrac{j}{2}$.  

Many integrals $\int_a^b V \mathcal{B}_j \mathcal{B}_k dt$ must be computed 
numerically (by a numerical integrator, such as "NIntegrate" in Mathematica), 
but for each $n \in \Z^+$ ($n \geq 3$) more than half of the 
entries of the $n \times n$ matrix $(\int_a^b V \mathcal{B}_j \mathcal{B}_k 
dt)_{j,k=1}^n$ can be determined to be zero simply by using the symmetry 
properties of the functions $V$ and $\mathcal{B}_j$.  For example, $V(a+t)=
V(a-t)$ and $\mathcal{B}_2(a+t)=\mathcal{B}_2(a-t)$, but $\mathcal{B}_3(a+t)=
-\mathcal{B}_3(a-t)$, and so $\int_a^b V \mathcal{B}_2 \mathcal{B}_3 dt$ must 
be zero.  

Let us list the eigenvalues of $\mathcal{L}_0$ in increasing order as 
\[ \lambda_0 < \lambda_1=-1 < \lambda_2=0 < \lambda_3 
\leq \lambda_4 \leq ... \to + \infty \; . \]
Each eigenvalue appears in this list the same number of times as the dimension 
of its eigenspace.  As noted in Proposition 
\ref{eigvalproperties}, $\lim_{j \to \infty} \lambda_j = + \infty$.  

\begin{Rmk}
In fact, all eigenvalues have multiplicity at most $2$.  No eigenspace can 
contain three independent eigenfunctions, as $\mathcal{L}_0$ is a second-order 
linear ODE.  
\end{Rmk}

For each $n \in \Z^+$, 
the $n \times n$ matrix $(\alpha_{jk})_{j,k=1}^n$ is symmetric, so it also 
has real eigenvalues, which we list in increasing order as 
\[ \lambda_0^{(n)} \leq \lambda_1^{(n)} \leq \lambda_2^{(n)} \leq ... \leq 
\lambda_{n-1}^{(n)} \; . \]
Rayleigh quotient characterizations for $\lambda_j$ and $\lambda_j^{(n)}$ prove that 
\begin{equation}\label{goingdown} \lambda_j^{(j+1)} \geq 
\lambda_j^{(j+2)} \geq \lambda_j^{(j+3)} \geq ... 
\geq \lambda_j \; . \end{equation} 
Thus the limit $\lim_{n \to \infty} \lambda_j^{(n)}$ exists and is greater 
than or equal to $\lambda_j$.  Using that 
$\{ \mathcal{B}_j=\mathcal{B}_j(t) \}_{j = 1}^\infty$ is an orthonormal basis 
of $\hat{\mathcal{F}}$, and so has dense span in $\hat{\mathcal{F}}$ with respect to 
the $L^2$ norm, further arguments 
with Rayleigh quotient characterizations give that in fact the limit is exactly 
equal to $\lambda_j$:  

\begin{theorem}\label{fromgeomdedicata} (\cite{r1}) 
$\lim_{n \to \infty} \lambda_j^{(n)} = \lambda_j$.  
\end{theorem}

Theorem \ref{fromgeomdedicata} is proven in \cite{r1}.  The essential facts behind its 
proof are well established and can be found in many sources 
(\cite{be}, \cite{chavel}, \cite{u} for example), but we reference \cite{r1} because 
the result is given there in a situation exactly analogous to the one here.  
This result is a simple variant of the basic theory in finite 
element methods (see the introductory chapters of 
\cite{strangfix} and \cite{brenscot}, for example), and it is essentially only a 
variant of the Ritz-Galerkin method (see \cite{mikhlin}, for example).  

This theorem now provides us with our second numerical method for 
estimating $\lambda_0$ 
(and any other eigenvalue of $\mathcal{L}_0$) by simply finding the smallest 
eigenvalue (and other eigenvalues) of the matrix 
$(\alpha_{jk})_{j,k=1}^n$ for sufficiently large $n$.  By the inequalities 
\eqref{goingdown}, it is clear that the estimates will be from above.  Numerical 
results are shown in Table \ref{thetable}, and they again confirm the primary numerical 
result in this paper.  

\begin{Rmk} Symmetries of the first eigenfunction can allow us to remove 
some $\mathcal{B}_j$, if we are only looking for the first eigenvalue $\lambda_0$.  
If $u$ is the eigenfunction for the first eigenvalue, then the symmetries of $u$ 
in Lemma \ref{Vsymmetry} imply that $u = \sum_{j \geq 1} a_j \mathcal{B}_j$ with 
$a_j=0$ when $j \geq 2$ is not an integer multiple of four.  Then we can 
more quickly estimate $\lambda_0$ by using subspaces of $\hat{\mathcal F}$ spanned 
by only $\mathcal{B}_1$ and $\mathcal{B}_{4k}$ for $k \in \Z^+$.  However, this 
shortcut will not work for estimating the other eigenvalues of $\mathcal{L}_0$.  
\end{Rmk} 

\begin{table}[htb]
    \begin{center}
    \begin{tabular}{l|c|c|c|c|c|c|c|c|c}
      & & Lemma & Lemma & First & Second & Second & Second & Second & Second \\
      & $\tfrac{b}{3}=$ & \ref{simplelemma1}'s lower & 
      \ref{simplelemma2}'s upper & method's & method's & 
      meth.'s & meth.'s & meth.'s & meth.'s
       \\
     $m$ &$=-a$ & bound & bound & estimate & estimate 
      & est. & est. & est. for & est. for \\
  &  & for $\lambda_0$ & for $\lambda_0$ & for $\lambda_0$ 
 & for $\lambda_0$ & for $\lambda_1$ & for $\lambda_2$ & 
  $\lambda_3$, $\lambda_4$ & $\lambda_5$, $\lambda_6$ \\ \hline
      $-1/4$ & $2.0137$ & $-2.25$ &$-1.0522$& $-1.25$ & $-1.245$ &$-0.992$&$0.00543$  &
           $1.44$        &$4.45$ \\ \hline
      $-1/2$ & $1.656$  & $-2.5$&$-1.3643$& $-1.5$ & $-1.4994$  &$-0.9987$ &$0.00102$& 
            $2.279$       &$6.75$ \\ \hline
      $-1$ & $1.3108$  & $-3$&$-1.9137$& $-2$ &  $-1.999$ &$-0.9993$ &$0.00062$& 
           $3.86$        &$11.021$ \\ \hline
      $-2$ & $1.0$  & $-4$ &$-2.9483$& $-3$ & $-2.999$ & $-0.9932$&$0.00615$& 
           $6.94$        &$19.26$ \\ \hline
      $-3$ & $0.8428$  & $-5$ &$-3.964$& $-4$ & $-3.999$ &$-0.999$ &$0.00069$& 
           $9.943$        &$27.3$ \\ \hline 
      $-10$ & $0.4849$  & $-12$&$-10.988$& $-11$ &  $-10.999$  &$-0.9974$ &$0.00256$& 
           $30.99$        &$83.46$ \\ \hline 
      $-20$ & $0.34683$  &$-22$ &$-20.994$& $-12$ &  $-20.999$  &$-0.983$ &$0.0168$& 
             $61.056$      & $163.61$
    \end{tabular}
    \end{center}
\caption{Estimates for $\lambda_0$ and other eigenvalues of 
$\mathcal{L}_0$ (computed using the values $H=1$ and $n=13$).}
\label{thetable}
\end{table}

\par
\setlength{\parindent}{0in}
\footnotesize
\vspace{0.2cm}
{\sc Wayne Rossman}:\, Department of Mathematics, Faculty of Science, 
Kobe University, 
Rokko, Kobe 657-8501, Japan. \, {\it E-mail}: 
wayne@math.kobe-u.ac.jp \, {\it web page}: 
www.math.kobe-u.ac.jp \par
\par

\end{document}